\documentclass[preprint,showpacs,showkeys,aps,amsmath,amssymb]{revtex4}%
\usepackage{amsthm}
\numberwithin{equation}{section}
\newtheorem{theorem}{Theorem}[section]
\newtheorem{corollary}[theorem]{Corollary}

\newtheorem{lemma}[theorem]{Lemma}
\newtheorem{proposition}[theorem]{Proposition}
\theoremstyle{definition}
\newtheorem{definition}[theorem]{Definition}
\newtheorem{example}[theorem]{Example}
\newtheorem{remark}[theorem]{Remark}
\renewcommand{\theenumi}{\roman{enumi}}

\begin{document}
\title[Existence of dynamics in quantum probability]{The existence
problem for dynamics of dissipative systems in quantum probability}
\author{Palle E. T. Jorgensen}
\email{jorgen@math.uiowa.edu}
\homepage{http://www.math.uiowa.edu/~jorgen/}
\thanks{Supported in part by the National Science Foundation under grants DMS-9987777 and DMS-0139473 (FRG).}
\affiliation{Department of Mathematics,
The University of Iowa,
14 MacLean Hall,
Iowa City, IA 52242-1419
U.S.A.}
  
\date{\today} 
  
\begin{abstract}
Motivated by existence problems for dissipative systems arising naturally
in lattice models from quantum statistical mechanics, we consider the
following $C^{\ast}$-algebraic setting:
A given hermitian dissipative mapping $\delta$ is densely defined
in a unital $C^{\ast}$-algebra $\mathfrak{A}$. The identity element in
$\mathfrak{A}$ is also in the domain of $\delta$. Completely dissipative maps
$\delta$ are defined by the requirement that the induced maps, $(a_{ij})
\rightarrow(\delta(a_{ij}))$, are dissipative on the $n$ by $n$
complex matrices over $\mathfrak{A}$ for all $n$. We establish the
existence of different types of maximal extensions of
completely dissipative maps. If the enveloping von Neumann algebra of
$\mathfrak{A}$ is injective, we show the existence of an
extension of $\delta$ which is the infinitesimal generator of a quantum
dynamical semigroup of completely positive maps in the von Neumann algebra.
If $\delta$ is a given
well-behaved $\ast$-derivation, then we show that each of the maps $\pm\delta$
is completely dissipative.

\end{abstract}


\maketitle

\section{Introduction\label{Int}}

Recent applications of the operator-theoretic approach to dissipative quantum
systems include \cite{FMRR00} and \cite{SW00}. For a more systematic approach,
see \cite{29}. Suppose we are given a one-parameter group of automorphisms
$\alpha_{t}\colon a\mapsto e^{itH}ae^{-itH}$ which acts on some set of
observables $a$, specified as a dense \textquotedblleft
local\textquotedblright\ subalgebra of a completed $C^{\ast}$-algebra. If we
then differentiate at $t=0$, we get the derivation $\delta\colon a\mapsto
i\left[  H,a\right]  =i\left(  Ha-aH\right)  $ which takes the form of a
formal commutator. The issue is complicated by the fact that the Hamiltonian
$H$ is typically an unbounded operator in statistical models, say infinite
lattice spin systems. In applications, it is $H$ that is given, and the
process must be run in reverse. By analogy to boundary value problems from
partial differential equations, we then expect to encounter an existence
problem for reconstructing the dynamics of the system from knowing only a
formula for $H$.

We adopt the $C^{\ast}$-$W^{\ast}$-formalism for the dynamics of infinite
quantum systems \cite{4,10,14,17,18,20,25,27}. For the special case of quantum
spin systems it is believed that the dynamics in the time reversible case is
given by an unbounded derivation of a suitable algebra $\mathfrak{A}$ of
observables \cite{25}. Depending on the range of the interaction, and the
number of dimensions of the spin lattice, it is possible to exponentiate the
infinitesimal derivation to a one-parameter group of automorphisms $\alpha
_{t}$ ($-\infty<t<\infty$) of $\mathfrak{A}$, or of the enveloping $W^{\ast}%
$-algebra $\mathfrak{A}^{\prime\prime}$ (see \cite{28}), or the $W^{\ast}%
$-algebra generated by a given invariant state \cite{4,13,21,23,26,27}.

It is known that (open) irreversible systems may be obtained as restrictions
of time-reversible systems, and it follows \cite{17} that the dynamics of the
open system is given mathematically by a semigroup $\tau_{t}$ ($0\leq
t<\infty$) of completely positive mappings of the $C^{\ast}$-algebra
$\mathfrak{A}$, or $W^{\ast}$-algebra $\mathfrak{A}^{\prime\prime}$. The
corresponding infinitesimal generator is completely dissipative. Completely
positive semigroups also play a role in quantum computing algorithms
\cite{Lin00}. The philosophy is that noise in the quantum processes dictates
the dissipative systems, as opposed to the conservative ones (which are
governed by one-parameter groups of automorphisms).

But in high lattice-dimensions, or for long-range interaction, there are
difficulties in exponentiating the infinitesimal generators. The determination
of the time-evolutions $\alpha_{t}$ (resp., $\tau_{t}$) seems to require
``extra boundary conditions'' \cite{3,4,21,26,27}. It is therefore a
meaningful foundational question, for a given completely dissipative
infinitesimal transformation $\delta$ in a $C^{\ast}$-algebra $\mathfrak{A} $,
to ask if it is always possible to extend $\delta$ to a transformation
$\tilde{\delta}$ which is the infinitesimal generator for a quantum dynamical
semigroup. Under the assumption that $\delta$ is hermitian, and that the
$W^{\ast}$-algebra $\mathfrak{A} ^{\prime\prime}$ is injective, we establish
the existence of a generator extension $\tilde{\delta}$. Our extension is thus
an algebraic parallel to Friederichs's extension for semibounded operators in
Hilbert space, or an analogue of Phillips's \cite{24} maximal dissipative
extension of the general dissipative operator in Hilbert space.

In earlier articles \cite{3,4,22,25} the uniqueness problem was considered for
the generator extension, $\delta\subset\tilde{\delta}$. But, just as is the
case for operators in Hilbert space (Friedrichs, Phillips), the extension is
generally not unique, reflecting the possibility of different ``boundary
conditions'' at infinity.

We refer the reader to the books \cite{10}, \cite{14} and \cite{27} for
details on the mathematical foundations of algebraic quantum theory.

The issues centering around the existence problem for the dynamical
one-parame\-ter groups, or semigroups, of quantum statistical mechanics are
perhaps best known in the setup of \emph{quantum spin systems}, as they are
treated in \cite{BrRoII}, \cite{JaPi01} and \cite{Rue01}.

\begin{example}
\label{ExaInt.1}The mathematical framework is rather general such as to allow
a wide variety of applications, including recent ones to nonequilibrium
statistical mechanics \cite{Rue01}. A countably infinite set $L$ (say a
lattice; it may be $\mathbb{Z}^{\nu}$ where $\nu$ is the lattice rank, or
dimension) is specified at the outset. Points $s\in L$ are sites at which
quantum spins are located. For each $s\in L$, let $\mathcal{H}_{s}$ be a
finite-dimensional complex Hilbert space, i.e., the spin vectors at site $s$;
and for a finite subset $\Lambda\subset L$, set%
\[
\mathcal{H}_{\Lambda}:=\bigotimes_{s\in\Lambda}\mathcal{H}_{s}.
\]
Then let $\mathfrak{A}_{\Lambda}$ be the $\ast$-algebra of all (bounded)
operators on $\mathcal{H}_{\Lambda}$. With the natural embedding%
\[
\mathfrak{A}_{\Lambda_{1}}\subset\mathfrak{A}_{\Lambda_{2}}\text{\qquad for
}\Lambda_{1}\subset\Lambda_{2}%
\]
given by%
\[
\mathfrak{A}_{\Lambda_{1}}\longmapsto\mathfrak{A}_{\Lambda_{1}}\otimes
1_{\Lambda_{2}\setminus\Lambda_{1}}\subset\mathfrak{A}_{\Lambda_{2}},
\]
we get the usual inductive limit $C^{\ast}$-algebra $\lim_{\Lambda
}\mathfrak{A}_{\Lambda}=:\mathfrak{A}$. A function $\Lambda\mapsto\Phi\left(
\Lambda\right)  =\Phi\left(  \Lambda\right)  ^{\ast}\in\mathfrak{A}_{\Lambda}$
defined on the finite subsets $\Lambda$ of $L$ is called an \emph{interaction}%
, and%
\begin{equation}
H_{\Phi}\left(  \Lambda\right)  =\sum_{X\subset\Lambda}\Phi\left(  X\right)
\label{eqInt.1}%
\end{equation}
is the associated local \emph{Hamiltonian}, where in (\ref{eqInt.1}), the
summation is over all finite subsets $X$ of $\Lambda$. Since $\mathfrak{A}%
_{\Lambda_{1}}$ and $\mathfrak{A}_{\Lambda_{2}}$ commute when $\Lambda_{1}%
\cap\Lambda_{2}=\varnothing$, it follows that%
\begin{equation}
\delta\left(  a\right)  =\lim_{\Lambda}\left[  H\left(  \Lambda\right)
,a\right]  \label{eqInt.2}%
\end{equation}
is well defined for all local observables $a$ in the dense $\ast$-subalgebra%
\[
\mathfrak{A}_{0}=\bigcup_{\Lambda\operatorname*{fin}}\mathfrak{A}_{\Lambda
}\text{\qquad in }\mathfrak{A}%
\]
where $\left[  \,\cdot\,,\,\cdot\,\right]  $ in (\ref{eqInt.2}) denotes the
usual commutator $\left[  b,a\right]  :=ba-ab$. Ruelle proved that, if $\Phi$
is translationally invariant, and if, for some $\lambda>0$,%
\begin{equation}
\sum_{n=0}^{\infty}e^{n\lambda}\sup_{s\in L}\sum_{\substack{s\in
X\operatorname*{fin}\\\operatorname*{card}X=n+1}}\left\Vert \Phi\left(
X\right)  \right\Vert <\infty,\label{eqInt.3}%
\end{equation}
then the $\ast$-derivation $\delta$ defined in (\ref{eqInt.2}) is the
infinitesimal generator of a one-parameter subgroup of $\ast$-automorphisms
$\left\{  \alpha_{t}\right\}  _{t\in\mathbb{R}}\subset\operatorname*{Aut}%
\left(  \mathfrak{A}\right)  $, which then satisfies%
\begin{equation}
\alpha_{t}\left(  a\right)  =\lim_{\Lambda\nearrow L}e^{itH\left(
\Lambda\right)  }ae^{-itH\left(  \Lambda\right)  }\label{eqInt.4}%
\end{equation}
for all $a\in\mathfrak{A}$ and $t\in\mathbb{R}$, i.e., it is approximately
inner. This means that, if $a\in\mathfrak{A}_{0}$, then%
\begin{equation}
\lim_{\substack{t\rightarrow0\\t\neq0}}t^{-1}\left(  \alpha_{t}\left(
a\right)  -a\right)  =\delta\left(  a\right)  .\label{eqInt.5}%
\end{equation}
Moreover, $\delta$ is, when extended from $\mathfrak{A}_{0}$, a closed $\ast
$-derivation, in the sense that the graph of $\delta$ is closed in
$\mathfrak{A}\times\mathfrak{A}$. But if $\Phi$ is not translationally
invariant, or if (\ref{eqInt.3}) is not known to hold, then no such conclusion
is within reach, and the issue of extensions of $\delta$ arises. We then ask
if some extension $\tilde{\delta}$ of $\delta$ to a generator of a
one-parameter group of automorphisms, or a semigroup of dissipations (see
details below), exists.
\end{example}

\section{Definitions and terminology\label{Def}}

Let $X$ and $Y$ be Banach spaces. Then the space of bounded linear operators
from $X$ to $Y$ is denoted $L(X,Y)$. The \emph{conjugate} (i.e.,
\emph{dual\/}) Banach space to $X$ is $L(X,\mathbb{C})$, and is denoted
$X^{\prime}$. If $\mathcal{H}$ is a Hilbert space, the $C^{\ast}$-algebra of
all bounded operators on $\mathcal{H}$ is denoted $B(\mathcal{H})$. Let
$\mathcal{L}$ be a linear subspace of $B(\mathcal{H})$ which is self-adjoint
and contains the identity operator $I$. With the order inherited from
$B(\mathcal{H})$, the subspace $\mathcal{L}$ gets the structure of an
\emph{operator system}, in the terminology of Effros \cite{7}. The full matrix
algebra $M_{n}$ of all complex $n$-by-$n$ matrices is also an operator system,
and so is $\mathcal{L}_{n}=\mathcal{L}\otimes M_{n}$. The elements in
$\mathcal{L}_{n}$ may be realized as $n$-by-$n$ matrices with entries from
$\mathcal{L}$, $(a_{ij})_{i,j=1}^{n}$, $a_{ij}\in\mathcal{L}$. If
$\mathcal{L}$ and $\mathcal{R}$ are operator systems and $\varphi
\colon\mathcal{L}\rightarrow\mathcal{R}$ is a linear mapping, then the induced
map $(a_{ij})\rightarrow(\varphi(a_{ij}))$ of $\mathcal{L}_{n}$ into
$\mathcal{R}_{n}$ is denoted $\varphi_{n}$. It is, in fact, $\varphi
\otimes\operatorname*{id}_{n}$. We say \cite{2} that $\varphi$ is
\emph{completely positive} (resp., \emph{completely contractive\/}) if
$\varphi_{n}$ is positive (resp., contractive) for all $n$. We say that
$\mathcal{R}$ is \emph{injective} if for every pair of operator systems,
$\mathcal{L}\subset\mathcal{L}_{1}$, and every completely positive map
$\varphi\colon\mathcal{L}\rightarrow\mathcal{R}$, there is a completely
positive extension $\psi\colon\mathcal{L}_{1}\rightarrow\mathcal{R}$. That is,
$\psi(x)=\varphi(x)$ for all $x\in\mathcal{L}$. If $\mathcal{R}$ is a von
Neumann algebra in a Hilbert space $\mathcal{H}$, it is known \cite{30,9} that
$\mathcal{R}$ is injective iff there is a norm-one projection of
$B(\mathcal{H})$ onto $\mathcal{R}$.

If $\mathfrak{A} $ is a $C^{\ast}$-algebra, it is known \cite{8} that
$\mathfrak{A} $ is nuclear iff the double conjugate (dual) $\mathfrak{A}
^{\prime\prime}$ is injective as a $W^{\ast}$-algebra. Connes showed \cite{9}
that a factor $\mathcal{R}$ on a separable Hilbert space is injective iff it
is matricial.

\section{Dissipative transformations\label{Dis}}

An operator $\delta$ in a Banach space $X$ is said to be dissipative \cite{23}
if one of the following three equivalent conditions is satisfied:

\begin{enumerate}
\item \label{Dis(1)}For all $x$ in the domain $\mathcal{D}(\delta)$ of
$\delta$, there is an element $f\in X^{\prime}$, depending on $x$, such that
$\| f\| =1$, $f(x)=\| x\| $, and $\operatorname{Re}f(\delta(x))\leq0$.

\item \label{Dis(2)}For all $x$ in $\mathcal{D}(\delta)$, and all $f\in
X^{\prime}$ satisfying $\| f\| =1$, and $f(x)=\| x\| $, the inequality
$\operatorname{Re}f(\delta(x))\leq0$ is valid.

\item \label{Dis(3)}For all $x$ in $\mathcal{D}(\delta)$, and all $\alpha
\in\mathbb{R}_{+}$, the inequality $\| x-\alpha\delta(x)\| \geq\| x\| $ holds.
\end{enumerate}

The proof of the equivalence can be found, for example, in \cite{3}, but the
equivalence can also be shown to be a consequence of the approximation idea in
Section \ref{Con} and Proposition \ref{Proposition10} in the present paper.

If $X$ is an operator system, we say that $\delta$ is \emph{completely
dissipative} if the induced mapping $\delta_{n}$ in $X_{n}$ is dissipative for
all $n=1,2,\dots$. Recall that $X_{n}=X\otimes M_{n}$, and $\delta_{n}%
\colon(x_{ij})\rightarrow(\delta(x_{ij}))$, with domain $\mathcal{D}%
(\delta_{n})=\{(x_{ij})\in X_{n}:x_{ij}\in\mathcal{D}(\delta)\}$.

Finally we say that the transformation $\delta$ is \emph{hermitian} if the
domain $\mathcal{D}(\delta)$, in the operator system $X$, is invariant under
the $\ast$-involution of $X$, and if $\delta(x^{\ast})=\delta(x)^{\ast}$ for
all $x\in\mathcal{D}(\delta)$.

If $\delta\colon X\rightarrow Y$ is merely a linear transformation between
Banach spaces $X$ and $Y$, with domain $\mathcal{D}(\delta)$ dense in $X$,
then the transposed (or conjugate) transformation $\delta^{\prime}$ is well
defined as a linear transformation $\delta^{\prime}\colon Y^{\prime
}\rightarrow X^{\prime}$ with domain $\mathcal{D}(\delta^{\prime})=\{f\in
Y^{\prime}:\exists\,g\in X^{\prime}$ s.t.\ $f(\delta(x))=g(x)$ for all
$x\in\mathcal{D}(\delta)\}$. For $f\in\mathcal{D}(\delta^{\prime})$,
$\delta^{\prime}(f)=g$. The domain $\mathcal{D}(\delta^{\prime})$ is
weak*-dense in $Y^{\prime}$ iff $\delta$ is closable. It is known \cite{23}
that dissipative operators are closable.

\section{Completely positive semigroups (Quantum dynamical
semigroups)\label{Com}}

Let $M$ be a $W^{\ast}$-algebra with predual $M_{\ast}$. Let $\tau_{t}$ be a
family of completely positive mappings of $M$ into itself, indexed by the time
parameter $t\in\lbrack0,\infty)$. Assume that $\tau_{0}$ is the identity
transformation in $M$, and that $\tau_{t}(\openone      )=\openone        $
for all $t\in\lbrack0,\infty)$, where $\openone        $ denotes the unit
element of the $W^{\ast}$-algebra $M$ in question. We assume further that the
semigroup law holds, $\tau_{t_{1}+t_{2}}=\tau_{t_{1}}\circ\tau_{t_{2}}$ for
$t_{1},t_{2}\in\lbrack0,\infty)$, and finally that each $\tau_{t}$ is a normal
mapping in $M$. Recall that normality is equivalent to the requirement that
the conjugate semigroup $\tau_{t}^{\prime}$ \cite{12} of $M^{\prime}$ leaves
invariant the subspace $M_{\ast}$. Finally we require continuity of each
scalar function, $t\rightarrow\varphi(\tau_{t}(a))$, for all $\varphi\in
M_{\ast}$ and $a\in M$. A semigroup which satisfies all the requirements above
is called a \emph{completely positive semigroup}. Because of the relevance to
quantum dynamics, we shall also call it a quantum dynamical semigroup
\cite{18}.

The \emph{infinitesimal generator} of a given completely positive semigroup
$(\tau_{t},M)$ is a, generally unbounded, transformation, denoted by $\zeta$,
in $M$. The domain of the generator $\zeta$ is given by
\[
\mathcal{D}(\zeta)=\{a\in M:\exists\,b\in M\text{ s.t.\ for all }t,\;\tau
_{t}(a)-a=%
{\textstyle\int_{0}^{t}}
\tau_{s}(b)\,ds\}.
\]

By definition $\zeta(a)=b$. It is easy to see \cite{12} that $\zeta
(a)=\frac{d\,}{dt}\tau_{t}(a)|_{t=0}$, where the derivative is taken in the
$\sigma(M,M^{\ast})$-topology. Finally note that infinitesimal generators are
completely dissipative.

\begin{example}
\label{ExaCom.1}It is known that the generator $\delta$ of a completely
positive semigroup $\left\{  \tau_{t}\right\}  _{t\in\mathbb{R}_{+}}$ on a
$C^{\ast}$-algebra $\mathfrak{A}$ is completely dissipative on a dense
subspace $\mathcal{D}$ in $\mathfrak{A}$; see \cite{Arv02a}. The following is
a \textquotedblleft canonical\textquotedblright\ example of this: it is built
on the $C^{\ast}$-algebra over the canonical commutation relations (CCR); see
\cite{Bha01}. Specifically, let $\mathcal{H}$ be a complex Hilbert space. Then
there is a $C^{\ast}$-algebra $\mathfrak{A}=\mathfrak{A}\left(  \mathcal{H}%
\right)  $ which is generated by the identity element $\openone$ and a family
of unitary elements $\left\{  u_{\xi}\mid\xi\in\mathcal{H}\setminus\left\{
0\right\}  \right\}  $ such that
\[
u_{\xi}u_{\eta}=e_{\mathstrut}^{\frac{i}{2}\operatorname{Im}
\left\langle\,\xi\mid\eta\,\right\rangle}%
u_{\xi+\eta}%
\]
for all $\xi,\eta\in\mathcal{H}$, with the understanding that $u_{0}%
=\openone$. Then it follows that there is a unique, completely positive
semigroup $\left\{  \tau_{t}\right\}  _{t\in\mathbb{R}_{+}}$ in $\mathfrak{A}%
$, such that%
\[
\tau_{t}\left(  u_{\xi}\right)  =e_{\mathstrut}^{-t\left\Vert \xi\right\Vert
_{\mathcal{H}}^{2}}u_{\xi}\text{\qquad for }\xi\in\mathcal{H}.
\]
Hence the subalgebra $\mathcal{D}\subset\mathfrak{A}$ spanned by the elements
$\left\{  u_{\xi}\mid\xi\in\mathcal{H}\right\}  $ is contained in the domain
of the generator $\delta$, and%
\begin{equation}
\delta\left(  u_{\xi}\right)  =-\left\Vert \xi\right\Vert _{\mathcal{H}}%
^{2}u_{\xi}.\label{eqgenerator.star}%
\end{equation}
It follows from the observation in \cite{Arv97} and \cite{Arv02b} that this
$\delta$ is completely dissipative with dense domain $\mathcal{D}$ in the
$C^{\ast}$-algebra $\mathfrak{A}$. That is, $\delta$ defined by%
\[
\delta\left(  a\right)  =\lim_{t\rightarrow0_{\rlap{$\scriptscriptstyle +$}}%
}\phantom{\scriptscriptstyle +}t^{-1}\left(  \tau_{t}\left(  a\right)
-a\right)  \text{\qquad(norm limit)}%
\]
is well defined for $a=u_{\xi}\in\mathcal{D}$, and (\ref{eqgenerator.star}) holds.
\end{example}

We now turn to the general existence problem.

\begin{theorem}
\label{Theorem1}Let $\mathfrak{A} $ be a $C^{\ast}$-algebra with unit
$\openone        $, and let $\delta$ be a completely dissipative
transformation in $\mathfrak{A} $ with dense domain $\mathcal{D}(\delta)$.
Assume $\openone  \in\mathcal{D}(\delta)$, $\delta(\openone        )=0$, and
further that $\delta$ is hermitian. Moreover assume that the double conjugate
\textup{(}dual\/\textup{)} $\mathfrak{A} ^{\prime\prime}$ is an injective
$W^{\ast}$-algebra. Then $\delta$ has an extension $\tilde{\delta}$ to an
ultraweakly densely defined transformation in $\mathfrak{A} ^{\prime\prime}$
which is at the same time the infinitesimal generator of a completely positive
semigroup of normal unital transformations in $\mathfrak{A} ^{\prime\prime}$.
\end{theorem}

We have divided the proof of Theorem \ref{Theorem1} into two main sections:
one is concerned with the analysis of the family of extensions of the
\emph{partial resolvent operator} $(I-\delta)^{-1}$. This analysis leads to a
distinguished set of contractive, and maximal, extensions which is associated
with a set of extensions $\tilde{\delta}$ of $\delta$. But $\tilde{\delta}$
turns out to be an operator in the enveloping $W^{\ast}$-algebra of
$\mathfrak{A} $. The generation properties of $\tilde{\delta}$ are analyzed in
the second section.

\section{Extensions of $(I-\delta)^{-1}$\label{Ext}}

We may assume that $\delta$ is in fact a closed operator in $\mathfrak{A} $.
(If not, it would be possible to replace $\delta$ by the closure $\bar{\delta
}$, and $\bar{\delta}$ will have the properties which were listed for $\delta$.)

This means that the linear space $\mathcal{S}=\operatorname*{Ran}
(I-\delta)=\{x-\delta(x):x\in\mathcal{D}(\delta)\}$ is closed in $\mathfrak{A}
$. In view of the (hermitian) assumption on $\delta$ we note that
$\mathcal{S}$ is also selfadjoint, and that $\openone        \in\mathcal{S}$.
The operator $R\colon\mathcal{S}\rightarrow\mathfrak{A} $ defined by
$x-\delta(x)\rightarrow x$, and denoted by $(I-\delta)^{-1}$, is completely
positive \cite[Prop.~1.2.8]{2}. Clearly $R(\openone        )=\openone        $.

We now consider the double dual to $\mathfrak{A} $, denoted by $\mathfrak{A}
^{\prime\prime}$, as a $W^{\ast}$-algebra $M$, and make the appropriate
identification (via the universal $\ast$-representation for $\mathfrak{A} $)
such that $\mathfrak{A} $ is regarded as a $C^{\ast}$-subalgebra of
$\mathfrak{A} ^{\prime\prime}$, and the pre-dual of $\mathfrak{A}
^{\prime\prime}$ is identified with the dual $\mathfrak{A} ^{\prime}$ of
$\mathfrak{A} $. (The reader is referred to \cite[\S 1.17, p.~42]{28} for
details.) Since $M=\mathfrak{A} ^{\prime\prime}$ (with the Arens
multiplication) is injective as a $W^{\ast}$-algebra, by the assumption, it
follows that a completely positive extension mapping $E\colon M\rightarrow M$
exists. If we regard $\mathfrak{A} $ as a subalgebra of $M$ (as we shall),
then the extension property is given by the identity
\begin{equation}
R(s)=E(s)\text{\qquad for all }s\in\mathcal{S}. \label{eq5.1}%
\end{equation}
Note that $\mathcal{S}\subset\mathfrak{A} $, so that $\mathcal{S}$ becomes a
subspace of $M$ with the above mentioned identification.

The completely positive transformations of $M$ into itself will be denoted by
$CP(M)$, and the space $L(M)$ of completely bounded linear transformations in
$M$ gets an ordering arising from the cone $CP(M)$. Indeed, for $F\in L(M)$ we
define $E\leq F$ by the requirement that $F-E\in CP(M)$. Among all the
particular extensions $F$ of $R$, $F\in L(M)$, such that $E\leq F$, we choose
by Zorn a maximal element $F_{0}$. (For the basic facts on topologies on
$CP(M)$ which are needed, the reader is referred to \cite[Ch.~1]{2}.)

This extension $F_{0}$, described above, has the special property of being
$1$--$1$. We first consider the restriction of $F_{0}$ to the positive
elements in $M$, $M_{+}$, that is. More precisely, we have the implication:
\begin{equation}
x\in M_{+},\;F_{0}(x)=0\Longrightarrow x=0. \label{eq5.2}%
\end{equation}
Let $\eta\colon M\rightarrow M/\mathcal{S}$ be the canonical linear quotient
mapping, and consider the cone $\mathcal{C}$ in the normed quotient space
$\mathcal{E}=\mathcal{M}/\mathcal{S}$ given by $\mathcal{C}=\eta(M_{+})$.

If the element $x$ in (\ref{eq5.2}) belongs to $\mathcal{S}$, then the
conditions $R(x)=F_{0}(x)=0$ imply $x=0$, since $R=(I-\delta)^{-1}$. Hence, we
shall assume that $x$ is not in $\mathcal{S}$. This means that $\eta
(x)\in\mathcal{C}$ defines a one-dimensional subspace $\{k\eta(x):k\in
\mathbb{C}\}$ in $\mathcal{E}$, and the functional $f\colon k\eta
(x)\rightarrow k$ is nonzero and positive. By Krein's theorem \cite[Thm.~1,
Ch.~3, p.~157]{1} $f$ extends to a positive functional $\tilde{f}$ on
$\mathcal{E}$, and we may define
\begin{equation}
F_{1}(y)=F_{0}(y)+\tilde{f}(\eta(y))\openone \text{\qquad for }y\in M.
\label{eq5.3}%
\end{equation}
We claim that $F_{1}$ is one of the extensions considered in the Zorn-process
which was described above. But $F_{0}\leq F_{1}$, and $F_{0}\neq F_{1}$,
contradicting the maximality of $F_{0}$---and so, we must have $x=0$,
concluding the proof of (\ref{eq5.2}). (Note that in (\ref{eq5.3}), instead of
the identity element $\openone        $ on the right-hand side of the
equation, we could have used any nonzero element in $M_{+}$. The corresponding
$F_{1}$-transformation would properly majorize $F_{0}$, and have its range
contained in $M$, since the range of $F_{0}$ falls in $M$.)

Since $F_{0}$ is completely positive, we have, in particular, $F_{0}(x^{\ast
})=F_{0}(x)^{\ast}$. So, to establish the identity $N(F_{0})=\{x\in
M:F_{0}(x)=0\}=0$, it is enough to show that the hermitian part of $N(F_{0})$
is zero. Since we have already considered positive elements, it only remains
to consider $x=x^{\ast}\in N(F_{0})$ satisfying $x\notin\mathcal{S}$. Choose a
positive real number $k$ such that $x_{k}=x+k\openone        \in M_{+}$. We
then have $F_{0}(x_{k})=k$ and $x_{k}\notin\mathcal{S}$. It is possible,
therefore, by Krein's theorem, to choose a positive functional $\tilde{f}$ on
$\mathcal{E}=M/\mathcal{S}$ satisfying $\tilde{f}(\eta(x_{k}))=l>0$. Then
define $F_{2}(y)=F_{0}(y)+\tilde{f}(\eta(y))\openone        $ for $y\in M$. It
is a simple matter to check that $F_{2}$ is one of the Zorn-extensions.
Indeed, $F_{0}\leq F_{2}$ since $\tilde{f}$ is chosen positive. Finally
$F_{2}(x_{k})=F_{0}(x_{k})+l\openone        >F_{0}(x_{k})$. This contradiction
to the maximality of $F_{0}$ concludes the proof. Since $N(F_{0})=0$, the
inverse $F_{0}^{-1}$ is defined on $F_{0}(M)=\{F_{0}(x):x\in M\}$.

We proceed to show that $F_{0}(M)$ is in fact dense in the $\sigma
(M,\mathfrak{A} ^{\prime})$-topology of $M$: First note that the extension
property (\ref{eq5.1}) for $F_{0}$ translates into:
\begin{equation}
F_{0}(x-\delta(x))=x\text{\qquad for }x\in\mathcal{D}(\delta), \label{eq5.4}%
\end{equation}
and the corresponding transposed mappings in $\mathfrak{A} ^{\prime}$
therefore satisfy:
\begin{equation}
(I-\delta^{\prime})F_{0}^{\prime}=I\text{\qquad(the identity operator in
}\mathfrak{A} ^{\prime}\text{).} \label{eq5.5}%
\end{equation}
Hence $F_{0}^{\prime}$ is $1$--$1$, and the desired density of $F_{0}(M)$
follows from the bi-polar theorem applied to the $\mathfrak{A} ^{\prime}$--$M$
duality. Note that in fact every extension of $R$ has dense range, because
condition (\ref{eq5.5}) is satisfied for the most general such extension.

Since $F_{0}$ is an extension of $(I-\delta)^{-1}$ it is clear that
$\tilde{\delta}=I-F_{0}^{-1}$ is therefore an extension of $\delta$.

\section{Generation properties of $\tilde{\delta}$\label{Gen}}

\addtolength{\textheight}{0.5\baselineskip}The operator $\tilde{\delta}$ is
closed and densely defined in the $\sigma$-topology of $M$. But $(I-\tilde
{\delta})^{-1}=F_{0}$, so we also have $\Vert x-\tilde{\delta}(x)\Vert
\geq\Vert x\Vert$ for all $x\in\mathcal{D}(\tilde{\delta})$. We proceed to
show that in fact%
\begin{equation}
\Vert kx-\tilde{\delta}(x)\Vert\geq k\Vert x\Vert\label{eq6.1}%
\end{equation}
for all $k>0$ and $x\in\mathcal{D}(\tilde{\delta})$. Indeed, let $\Lambda$
denote the set of $k>0$ such that the inequality (\ref{eq6.1}) is satisfied
for all $x\in\mathcal{D}(\tilde{\delta})$. Then we have seen that $k=1$
belongs to $\Lambda$. It turns out that $\Lambda$ is both open and closed as a
subset of $\mathbb{R}_{+}$, and our result follows by connectedness.

To show openness, suppose first that $k_{0}\in\Lambda$, and that
$k\in\mathbb{R}_{+}$ satisfies $\left|  k-k_{0}\right|  <k_{0}$. We than use
(\ref{eq6.1}), for $k_{0}$, in estimating the terms in the Neumann expansion
for $(kI-\tilde{\delta})^{-1}$, taken around the point $k_{0}$. Due to the
assumption $\left|  k-k_{0}\right|  <k_{0}$, the Neumann series is convergent,
and does indeed define a bounded inverse $R(k,\tilde{\delta})$ to
$kI-\tilde{\delta}$. Termwise estimation gives $\| R(k,\tilde{\delta})\| \leq
k^{-1}$, and it follows that (\ref{eq6.1}) is satisfied in a neighborhood of
$k_{0}$.

Consider next a convergent sequence of points $k_{n}\rightarrow k_{0}$ with
$k_{n}\in\Lambda$ and $k_{0}\in\mathbb{R}_{+}$. By assumption the resolvent
operators $R(k_{n},\tilde{\delta})=(k_{n}I-\tilde{\delta})^{-1}$ exist, and
they therefore satisfy the resolvent identity:%
\[
R(k_{n},\tilde{\delta})-R(k_{m},\tilde{\delta})=(k_{n}-k_{m})R(k_{n}%
,\tilde{\delta})R(k_{m},\tilde{\delta}),
\]
as well as the estimate $\Vert R(k_{n},\tilde{\delta})\Vert\leq k_{n}^{-1}$.
It follows that the norm-limit $\tilde{R}=\lim_{n}R(k_{n},\tilde{\delta})\in
L(M)$ exists, and it is trivial to check that $\tilde{R}$ defines a bounded
inverse to $k_{0}I-\tilde{\delta}$. The estimate (\ref{eq6.1}) for $k_{0}$ is
now implied in the limit by $\Vert\tilde{R}\Vert\leq k_{0}^{-1}$. Hence
$\Lambda$ is closed, and the argument is completed.

We have shown that the operator $\tilde{\delta}$ in $M$ is dissipative and
closed in the $\sigma(M,\mathfrak{A} ^{\prime})$-topology. It is, of course,
also closed in the norm-topology, and it can be showr that $\mathcal{D}%
(\tilde{\delta})$ is norm-dense. It follows by semigroup theory \cite{19,23}
that $\tilde{\delta}$ is the infinitesimal generator of a strongly continuous
semigroup $\tau_{t}$ ($0\leq t<\infty$) of contraction operators in the Banach
space $M$.

To show that each $\tau_{t}$ is a normal transformation we consider the
adjoint semigroup $\tau_{t}^{\prime}$ (cf.\ \cite{12}) in the norm-dual
$M^{\prime}$ and show that $\tau_{t}^{\prime}$ leaves $\mathfrak{A} ^{\prime}$
invariant. Note that $\mathfrak{A} ^{\prime}$ is being identified with the
predual of the $W^{\ast}$-algebra $M$, so that we may regard it as a subspace
of $M^{\prime}$.

Let $\tilde{\delta}^{\prime}$ (resp., $F_{0}^{\prime}$) denote the transposed
operators to $\tilde{\delta}$ (resp., $F_{0}$) with respect to the
$M$--$M^{\prime}$ duality. It follows by operator theory that $\tilde{\delta
}^{\prime}$ is the generator of $\tau_{t}^{\prime}$, and that $(I-\tilde
{\delta}^{\prime})^{-1}=F_{0}^{\prime}$. From the construction of $F_{0}$ we
now deduce that $\mathfrak{A} ^{\prime}$ is invariant under $F_{0}^{\prime}$.
Indeed, recall that $\delta^{\prime}$ denotes the transposed transformation to
$\delta$ with respect to the $\mathfrak{A} $--$\mathfrak{A} ^{\prime}$
duality. By definition $\mathcal{D}(\delta^{\prime})=\{a^{\prime}%
\in\mathfrak{A} ^{\prime}:\exists\,b^{\prime}\in\mathfrak{A} ^{\prime
},\;\left\langle b^{\prime},x\right\rangle =\left\langle a^{\prime}%
,\delta(x)\right\rangle $ for all $x\in\mathcal{D}(\delta)\}$. But for
$a^{\prime}\in\mathfrak{A} ^{\prime}$ and $x\in\mathcal{D}(\delta)$ we have
$\left\langle F_{0}^{\prime}(a^{\prime}),x-\delta(x)\right\rangle
=\left\langle a^{\prime},x\right\rangle $. Hence, $F_{0}^{\prime}(a^{\prime
})\in\mathcal{D}(\delta^{\prime})\subset\mathfrak{A} ^{\prime}$ by
(\ref{eq5.5}).

An application of the Neumann expansion to $(I-\frac{t}{n}\tilde{\delta
}^{\prime})^{-1}$ shows that $\mathfrak{A}^{\prime}$ is also invariant under
this operator for all $t\geq0$, $n\in\mathbb{Z}_{+}$. But $\tau_{t}^{\prime}$
is obtained as a weak*-limit of these operators ($n\rightarrow\infty$), and
the desired invariance $\tau_{t}^{\prime}(\mathfrak{A} ^{\prime}%
)\subset\mathfrak{A} ^{\prime}$ follows.

A final application of the Neumann series, now to the operators $(I-\frac
{t}{n}\tilde{\delta})^{-1}$, shows that $\tau_{t}$ is completely positive in
$M$ for all $t\geq0$. Indeed $(I-\frac{t}{n}\tilde{\delta})^{-1}$ may be
expanded in a norm-convergent power series in the completely positive operator
$F_{0}=(I-\tilde{\delta})^{-1}$, and $\tau_{t}=\lim_{n\rightarrow\infty
}(I-\frac{t}{n}\tilde{\delta})^{-1}$.

\section{The inequality $\delta(x^{\ast}x)\geq\delta(x)^{\ast}x+x^{\ast}%
\delta(x)$\label{Ine}}

It was shown in \cite{16} that if $\delta$ is a bounded hermitian linear map
in a $C^{\ast}$-algebra $\mathfrak{A} $, then the following two conditions are
equivalent:%
\begin{equation}
e^{t\delta}(x^{\ast}x)\geq e^{t\delta}(x^{\ast})e^{t\delta}(x),\qquad
\forall\,x\in\mathfrak{A} ,\;t\in\mathbb{R}_{+}, \label{eq7.1}%
\end{equation}
and%
\begin{equation}
\delta(x^{\ast}x)\geq\delta(x^{\ast})x+x^{\ast}\delta(x),\qquad\forall
\,x\in\mathfrak{A} . \label{eq7.2}%
\end{equation}
For unbounded $\mathfrak{A} $ the situation is not as well understood. It is
therefore of interest to study the connection between the property
(\ref{eq7.2}) for $\delta$, and the other conditions which are customarily
used in the applications of unbounded dissipative mappings in operator
algebras to quantum dynamics.

\begin{theorem}
\label{Theorem2}Let $\mathfrak{A} $ be a $C^{\ast}$-algebra with unit
$\openone       $, and let $\delta$ be a completely dissipative transformation
in $\mathfrak{A} $ with dense domain $\mathcal{D}(\delta)$ . Assume $\openone
\in\mathcal{D}(\delta)$, and $\delta(\openone       )=0$.\renewcommand{\theenumi}{\alph{enumi}}

\begin{enumerate}
\item \label{Theorem2(1)}Let $x\in\mathcal{D}(\delta)$ and assume that
$x^{\ast}x\in\mathcal{D}(\delta)$. Then%
\begin{equation}
\delta(x^{\ast}x)\geq\delta(x)^{\ast}x+x^{\ast}\delta(x). \label{eq7.3}%
\end{equation}

\item \label{Theorem2(2)}Suppose both $x$ and $x^{\ast}$ belong to
$\mathcal{D}(\delta)$. Then $\delta(x^{\ast})=\delta(x)^{\ast}$.
\end{enumerate}
\end{theorem}

The following results are corollaries to the proofs of Theorems \ref{Theorem1}
and \ref{Theorem2}.

\begin{corollary}
\label{Corollary3}Let $\mathfrak{A} $ be a $C^{\ast}$-algebra with unit
$\openone       $, and let $\delta$ be completely dissipative in $\mathfrak{A}
$ with dense domain $\mathcal{D}(\delta)$, $\openone       \in\mathcal{D}%
(\delta)$, $\delta(\openone       )=0$. \renewcommand{\theenumi}{\alph{enumi}}

\begin{enumerate}
\item \label{Corollary3(1)}If $\mathfrak{A} \subset B(\mathcal{H})$ for some
Hilbert space $\mathcal{H}$, then there is a sequence of completely positive
maps $E_{n}\colon\mathfrak{A} \rightarrow B(\mathcal{H})$, $E_{n}(\openone
)=\openone       $, such that the following norm-convergence holds:
\begin{equation}
E_{n}(x)\longrightarrow x\text{\qquad for }x\in\mathfrak{A} , \tag{i}%
\end{equation}
and
\begin{equation}
n(E_{n}(x)-x)\longrightarrow\delta(x)\text{\qquad for }x\in\mathcal{D}%
(\delta). \tag{ii}%
\end{equation}

\item \label{Corollary3(2)}If $\mathcal{D}(\delta)$ is hermitian, then
$\delta$ is hermitian as well, i.e., $\delta(x^{\ast})=\delta(x)^{\ast}$ for
all $x\in\mathcal{D}(\delta)$, and it is then possible, for each $n$, to
choose $E_{n}$ to be $1$--$1$ with dense range.

\item \label{Corollary3(3)}Let $\delta$ and $\mathfrak{A} $ be as in
\textup{(\ref{Corollary3(1)})}, and let $\pi\colon\mathfrak{A} \rightarrow
B(\mathcal{K})$ be a representation of $\mathfrak{A} $ in a Hilbert space
$\mathcal{K}$. Then there exists a sequence $E_{n}\in CP(\mathfrak{A}
,B(\mathcal{K}))$ such that the following norm convergence holds:
\begin{equation}
E_{n}(x)\longrightarrow\pi(x)\text{\qquad for }x\in\mathfrak{A} ,
\tag{i$^\prime$}%
\end{equation}
and
\begin{equation}
n(E_{n}(x)-\pi(x))\longrightarrow\pi(\delta(x))\text{\qquad for }%
x\in\mathcal{D}(\delta). \tag{ii$^\prime$}%
\end{equation}

\end{enumerate}
\end{corollary}

\begin{proof}
[Proofs]We consider again the range subspace $\mathcal{S}=\operatorname*{Ran}
(I-\delta)=\{x-\delta(x):x\in\mathcal{D}(\delta)\}$. As in the proof of
Theorem \ref{Theorem1} note that $R=(1-\delta)^{-1}\colon\mathcal{S}%
\rightarrow\mathfrak{A} $ is completely contractive, and $R(\openone
)=\openone       $. If $\mathfrak{A} $ is considered as a subalgebra of
$B(\mathcal{H})$, where $\mathcal{H}$ is the Hilbert space of the universal
representation, then there is, by Arveson's extension theorem \cite[Theorem
1.2.9]{2} a completely positive mapping $E\colon\mathfrak{A} \rightarrow
B(\mathcal{H})$ such that
\begin{equation}
R(s)=E(s)\text{\qquad for all }s\in\mathcal{S}. \label{eq7.4}%
\end{equation}

If for each $n=1,2,\dots$ the operator $\delta$ is replaced by $n^{-1}\delta$,
then the above argument yields a completely positive map $E_{n}\colon
\mathfrak{A} \rightarrow B(\mathcal{H})$ such that $E_{n}$ is an extension of
the partially defined operator $(I-n^{-1}\delta)^{-1}$.

We claim that the sequence $(E_{n})$ satisfies conditions (i) and (ii) which
are listed in Corollary \ref{Corollary3}(\ref{Corollary3(1)}). Indeed, for $x$
in dense $\mathcal{D}(\delta)$ we have $E_{n}(x-n^{-1}\delta(x))=x$, and
therefore
\begin{equation}
E_{n}(x)=n^{-1}E_{n}(\delta(x))+x, \label{eq7.5}%
\end{equation}
and
\begin{equation}
E_{n}(\delta(x))=n(E_{n}(x)-x). \label{eq7.6}%
\end{equation}
Passing to the limit in (\ref{eq7.5}), we get (i) for the special case
$x\in\mathcal{D}(\delta)$, but then also for all $x$ in $\mathfrak{A} $ by a
$3$-$\varepsilon$ argument since each $E_{n}$ is contractive. The result (ii)
of Corollary \ref{Corollary3}(\ref{Corollary3(1)}) is now an immediate
consequence of (\ref{eq7.6}).

Returning to the proof of Theorem \ref{Theorem2}, we note that
(\ref{Theorem2(2)}) is trivial from (ii). Indeed, for $x$ and $x^{\ast}$ in
$\mathcal{D}(\delta)$ we have
\[
\delta(x^{\ast})=\lim n(E_{n}(x^{\ast})-x^{\ast})=\lim_{n}(n(E_{n}%
(x)-x))^{\ast}=\delta(x)^{\ast}.
\]
The proof of Theorem \ref{Theorem2}(\ref{Theorem2(1)}) is based on both (i)
and (ii), together with the Kadison-Schwarz inequality for $E_{n}$: Suppose
$x\in\mathcal{D}(\delta)$ and $x^{\ast}x\in\mathcal{D}(\delta)$. Then
$\delta(x^{\ast}x)=\lim n(E_{n}(x^{\ast}x)-x^{\ast}x)$. For each term on the
right-hand side we have:
\begin{align}
n(E_{n}(x^{\ast}x)-x^{\ast}x)  &  \geq n(E_{n}(x)^{\ast}E_{n}(x)-x^{\ast
}x)\label{eq7.7}\\
&  =\frac{1}{2}((n(E_{n}(x)-x))^{\ast}(E_{n}(x)+x)\nonumber\\
&  \qquad+(E_{n}(x)+x)^{\ast}n(E_{n}(x)-x))\nonumber\\
\longrightarrow &  \frac{1}{2}(\delta(x)^{\ast}(2x)+(2x)^{\ast}\delta
(x))=\delta(x)^{\ast}x+x^{\ast}\delta(x),\nonumber
\end{align}
where the last convergence $\longrightarrow$ is based on (i) and (ii) from
Corollary \ref{Corollary3}(\ref{Corollary3(1)}). Since $\delta(x^{\ast}x)$ is
obtained in the limit on the left, the desired inequality (\ref{eq7.3}) in
(\ref{Theorem2(1)}) of Theorem \ref{Theorem2} follows.

Only part (\ref{Corollary3(2)}) of the corollary remains. The technique from
the proof of Theorem \ref{Theorem1} is applied here. We go back to the
extension $E$ from (\ref{eq7.4}) in the beginning of the present proof.
Consider the ordering on all the extensions $F$ of $R$, $F\in L(\mathfrak{A}
,B(\mathcal{H}))$, which is induced by the cone $CP(\mathfrak{A}
,B(\mathcal{H}))$, and choose by Zorn a particular extension $F$, $E\leq F$,
which is maximal. The argument from the proof of Theorem \ref{Theorem1} then
shows that $F$ is $1$--$1$, and the range $\operatorname*{Ran}(F)$ is dense.
It follows that the operator $\tilde{\delta}=I-F^{-1}\colon\operatorname*{Ran}%
(F)\rightarrow\mathfrak{A} $ exists and satisfies $\tilde{\delta}%
(x)=\delta(x)$ for all $x\in\mathcal{D}(\delta)$.

If $\alpha$ is a positive real number, then the same construction may be
carried out for the transformation $\alpha\delta$, instead of $\delta$. Hence
we get completely positive unital maps $F_{\alpha}$ such that the inverse
$F_{\alpha}^{-1}$ exists for each $\alpha$, and the domain of $I-F_{\alpha
}^{-1}$ contains $\mathcal{D}(\delta)$. Moreover $\tilde{\delta}_{\alpha
}=I-F_{\alpha}^{-1}$ satisfies $\tilde{\delta}_{\alpha}(x)=\delta(x)$ for
$x\in\mathcal{D}(\delta)$. To get a sequence of mappings satisfying the
conditions in Corollary \ref{Corollary3}(\ref{Corollary3(2)}), we need only
take $E_{n}=F_{n^{-1}}$ in the special case $\alpha=n^{-1}$.

The proof of part (\ref{Corollary3(3)}) in the corollary is parallel to
(\ref{Corollary3(1)}) with the following modification: Arveson's extension
theorem is now applied to the mapping $\pi\circ(I-\delta)^{-1}\colon
\mathcal{S}\rightarrow B(\mathcal{K})$.
\end{proof}

\section{The implementation problem\label{Imp}}

The conclusion (ii$^{\prime}$) in Corollary \ref{Corollary3}%
(\ref{Corollary3(3)}) is of interest when one wants to implement the
transformation $\delta$ by a dissipative operator in Hilbert space. In
particular, one is interested in implementing a completely dissipative
$\delta$-operator by a dissipative Hilbert-space operator. We shall establish
a clear two-way connection between the dissipative notion for $\delta$, and
for the implementing Hilbert-space operator.

\begin{theorem}
\label{Theorem5}Let $\mathfrak{A}$ be a $C^{\ast}$-algebra with unit
$\openone     $, and let $\delta$ be a completely dissipative transformation
in $\mathfrak{A}$ with dense domain $\mathcal{D}(\delta)$. Assume $\openone
\in\mathcal{D}(\delta)$ and $\delta(\openone     )=0$. Let $\omega$ be a state
of $\mathfrak{A}$, and let $(\pi_{\omega},\mathcal{K}_{\omega},\Omega)$ be the
corresponding GNS representation of $\mathfrak{A}$. Let $\tilde{\omega}$ be
the vector state on $B(\mathcal{K}_{\omega})$ given by the cyclic vector
$\Omega$, i.e., $\tilde{\omega}(X)=\left\langle\, X\Omega\mid\Omega\,\right\rangle
$ for $X\in B(\mathcal{K}_{\omega})$, and assume that it is possible to choose
the sequence $(E_{n})\subset CP(\mathfrak{A},B(\mathcal{K}_{\omega}))$ from
Corollary \textup{\ref{Corollary3}(\ref{Corollary3(3)})} in such a manner
that
\begin{equation}
\tilde{\omega}(E_{n}(x))=\omega(x)\text{\qquad for all }x\in\mathfrak{A}.
\label{eq9.1}%
\end{equation}

Then there is a dissipative operator $L_{\omega}$ in $\mathcal{K}_{\omega}$
such that
\begin{equation}
\pi_{\omega}(\delta(x))\Omega=L_{\omega}(\pi_{\omega}(x)\Omega)\text{\qquad
for all }x\in\mathcal{D}(\delta). \label{eq9.2}%
\end{equation}

\end{theorem}

\begin{proof}
Let $\pi=\pi_{\omega}$, $\mathcal{K}=\mathcal{K}_{\omega}$, and let
$(E_{n})\subset CP(\mathfrak{A},B(\mathcal{K}))$ be a sequence which, along
with the conditions listed in Corollary \ref{Corollary3}(\ref{Corollary3(3)}),
also fulfills the invariance restriction (\ref{eq9.1}) of the present theorem.
For each $n$ define an operator $C_{n}$ in $\mathcal{K}$ as follows:%
\[
C_{n}(\pi(x)\Omega)=E_{n}(x)\Omega,\qquad x\in\mathfrak{A}.
\]
Then
\begin{align*}
\Vert C_{n}\pi(x)\Omega\Vert^{2}  &  =\Vert E_{n}(x)\Omega\Vert^{2}%
=\tilde{\omega}(E_{n}(x)^{\ast}E_{n}(x))\\
&  \leq\tilde{\omega}(E_{n}(x^{\ast}x))=\omega(x^{\ast}x)=\left\langle\,
\pi(x^{\ast}x)\Omega\mid\Omega\,\right\rangle =\Vert\pi(x)\Omega\Vert^{2},
\end{align*}
where the norm is that of $\mathcal{K}$, and where the Schwarz inequality is
applied to $E_{n}$. It follows that $C_{n}$ is well defined, and that it
extends by limits (in $\mathcal{K}$) to a contraction operator, $C_{n}\in
B(\mathcal{K})$, $\Vert C_{n}\Vert\leq1$.

By Corollary \ref{Corollary3}(\ref{Corollary3(3)})(ii$^{\prime}$), we then
have
\begin{align*}
\pi(\delta(x))\Omega &  =\lim n(E_{n}(x)\Omega-\pi(x)\Omega)\\
&  =\lim n(C_{n}(\pi(x)\Omega)-\pi(x)\Omega)\\
&  =\lim n(C_{n}-I)\pi(x)\Omega\text{\qquad for }x\in\mathcal{D}(\delta).
\end{align*}
As a consequence, the following quadratic form on $\mathcal{K}$:%
\[
\pi(x)\Omega,\;\pi(y)\Omega\longrightarrow\lim\left\langle\, n(C_{n}%
-I)\pi(x)\Omega\mid\pi(y)\Omega\,\right\rangle _{\mathcal{K}}%
\]
is well defined. Using the contractive property of $C_{n}$, it is easy to show
that this quadratic form is given by a dissipative operator $L$; that is to
say%
\[
\lim\left\langle\, n(C_{n}-I)\pi(x)\Omega\mid\pi(y)\Omega\,\right\rangle
=\left\langle\, L\pi(x)\Omega\mid\pi(y)\Omega\,\right\rangle .
\]
Since the limit on the left is also equal to the inner product
\[
\left\langle\,
\pi(\delta(x))\Omega\mid\pi(y)\Omega\,\right\rangle ,
\]
the identity
(\ref{eq9.2}) of the theorem follows.
\end{proof}

\section{A condition for complete dissipativeness\label{Con}}

In applications \cite{14,20,27} it is often possible to determine the
derivation $\delta$ in a particular representation. If moreover the derivation
is known to be implemented by a dissipative operator in the corresponding
Hilbert space, then it follows in special cases that $\delta$ itself is
completely dissipative.

\begin{theorem}
\label{Theorem6}Let $\mathfrak{A} $ be a $C^{\ast}$-algebra with unit
$\openone     $ and let $\delta$ be a densely defined transformation in
$\mathfrak{A} $ such that $\openone      \in\mathcal{D}(\delta)$ and
$\delta(\openone      )=0$. Let $\omega$ be a state on $\mathfrak{A} $ such
that $\delta$ is implemented by a dissipative Hilbert-space operator $L$ in
the representation $\pi_{\omega}$. Assume moreover that $\pi_{\omega}$ is
faithful, and that $L\Omega=0$, where $\Omega$ denotes the cyclic vector in
the GNS representation. Then $\delta$ is completely dissipative on its domain.
\end{theorem}

\begin{proof}
Let $\mathcal{H}=\mathcal{H}_{\omega}$ be the Hilbert space of the faithful
representation $\pi_{\omega}$ and let $L$ be the operator in $\mathcal{H}$
which is assumed to exist, satisfying conditions (\ref{Theorem6proof(1)}) and
(\ref{Theorem6proof(2)}) below:

\begin{enumerate}
\item \label{Theorem6proof(1)}The domain of $L$ is $\pi_{\omega}%
(\mathcal{D}(\delta))\Omega$, and $L$ is a dissipative operator in the Hilbert
space $\mathcal{H}$;

\item \label{Theorem6proof(2)}$L$ implements $\delta$ in the representation
$\pi_{\omega}$, which is equivalent to the requirement that $L^{\ast}$ is
defined on $\pi_{\omega}(\mathcal{D}(\delta))\Omega$, and that on this domain
the following operator identity is valid:%
\begin{equation}
\pi(\delta(a))=L\pi(a)+\pi(a)L^{\ast}\text{\qquad for all }a\in\mathcal{D}%
(\delta). \label{eq10.1}%
\end{equation}

\end{enumerate}

We show first that $\delta$ must necessarily be a dissipative operator.
Indeed, by Phillips's theorem \cite[Thm.~1.1.3]{24} an extension $\tilde{L}$
of $L$ exists which is the infinitesimal generator of a strongly continuous
semigroup $S(t)$ of contraction operators in the Hilbert space $\mathcal{H}$.
We note that $S(t)$ implements a semigroup $\sigma(t)$ of positive mappings in
$B(\mathcal{H})$, given by%
\begin{equation}
\sigma(t)(A)=S(t)AS(t)^{\ast} \label{eq10.2}%
\end{equation}
for all $t\in\lbrack0,\infty)$ and $A\in B(\mathcal{H})$. By semigroup theory
we note that the generator ($\zeta$ say) of $\sigma(t)$ is dissipative, so the
following estimate holds:%
\begin{equation}
\| A-\alpha\zeta(A)\| \geq\| A\| \label{eq10.3}%
\end{equation}
for all $\alpha\in\lbrack0,\infty)$ and $A\in\mathcal{D}(\zeta)$ .

If $\delta_{\omega}$ denotes the operator $\pi_{\omega}(a)\rightarrow
\pi_{\omega}(\delta(a))$ with domain $\pi_{\omega}(\mathcal{D}(\delta))$, then
we claim (easy proof) that%
\begin{equation}
\delta_{\omega}(A)=\zeta(A)\text{\qquad for all }A\in\mathcal{D}%
(\delta_{\omega}), \label{eq10.4}%
\end{equation}
and the known estimate (\ref{eq10.3}) above then implies%
\begin{equation}
\| \pi_{\omega}(a)-\alpha\pi_{\omega}(\delta(a))\| \geq\| \pi_{\omega}(a)\|
\label{eq10.5}%
\end{equation}
for $a\in\mathcal{D}(\delta)$ and $\alpha\in\lbrack0,\infty)$. But
$\pi_{\omega}$ is faithful (and hence isometric), so (\ref{eq10.5}) is in fact
equivalent to the dissipation estimate
\[
\| a-\alpha\delta(a)\| \geq\| a\|
\]
for the operator $\delta$ itself.

For each $n=1,2,\dots$, we now consider the tensor-product construction of the
$C^{\ast}$-algebra $\mathfrak{A} $ with the $n$-by-$n$ complex matrices
$M_{n}$; and we define $\mathfrak{A} _{n}=\mathfrak{A} \otimes M_{n}$,
$\delta_{n}=\delta\otimes\operatorname*{id}_{n}$, the operator obtained by
application of $\delta$ to each entry $a_{ij}$ in the matrix representation of
elements in $\mathfrak{A} _{n}$, $\omega_{n}=\omega\otimes\operatorname*{tr}%
_{n}$ where $\operatorname*{tr}_{n}$ denotes the normalized trace on $M_{n}$,
$\pi_{\omega_{n}}$: the GNS representation of $\mathfrak{A} _{n}$ associated
to $\omega_{n}$.

The problem is to show that each of the operators $\delta_{n}$ is dissipative.
We show that in fact $\delta_{n}$ is implemented by a dissipative
Hilbert-space operator in the representation $\pi_{\omega_{n}}$. Hence, the
first part of the proof applies and yields the conclusion of the claim since
each representation $\pi_{\omega_{n}}$ is faithful, being the tensor product
of faithful representations.

Let $\mathcal{H}_{n}$ denote the representation Hilbert space of $\pi
_{\omega_{n}}$. We proceed to find a dissipative operator $L_{n}$ in
$\mathcal{H}_{n}$ such that $\delta_{n}$ is implemented by $L_{n}$. In view of
(\ref{eq10.1}) this means that%
\[
\pi_{\omega_{n}}(\delta_{n}(a))=L_{n}\pi_{\omega_{n}}(a)+\pi_{\omega_{n}%
}(a)L_{n}^{\ast}%
\]
for all $a\in\mathcal{D}(\delta_{n})=\mathcal{D}(\delta_{n})\otimes M_{n}$
(algebraic tensor product) $\subset\mathfrak{A}_{n}$ as an operator identity
on $\pi_{\omega_{n}}(\mathcal{D}(\delta_{n}))\Omega_{n}\subset\mathcal{H}_{n}%
$. Here $\Omega_{n}$ denotes the cyclic vector for the representation
$\pi_{\omega_{n}}$, i.e.,
\begin{equation}
\omega_{n}(a)=\left\langle\, \pi_{\omega_{n}}(a)\Omega_{n}\mid\Omega
_{n}\,\right\rangle \text{\qquad for all }a\in\mathfrak{A}_{n}. \label{eq10.6}%
\end{equation}

Our next step is the verification of the following:%
\begin{align}
&  \operatorname{Re}\omega_{n}(a^{\ast}\delta_{n}(a))\leq0\text{\qquad for all
}a\in\mathcal{D}(\delta_{n}),\label{eq10.7}\\
&  L_{n}\Omega_{n}=0,\label{eq10.8}\\
&  \omega_{n}(a^{\ast}\delta_{n}(a))=\left\langle\, L_{n}\pi_{\omega_{n}%
}(a)\Omega_{n}\mid\pi_{\omega_{n}}(a)\Omega_{n}\,\right\rangle \text{\qquad for
}a\in\mathcal{D}(\delta_{n}). \label{eq10.9}%
\end{align}
It will follow from (\ref{eq10.7}) and (\ref{eq10.9}) that an implementing
operator $L_{n}$ satisfying (\ref{eq10.8}) must necessarily be dissipative.

Note that (\ref{eq10.8}) is verified for $n=1$ by assumption. Hence
$\omega(a^{\ast}\delta(a))=
\left\langle\, \pi(\delta(a))\Omega\mid
\pi(a)\Omega\,\right\rangle =\left\langle\, L\pi(a)\Omega+\pi(a)L^{\ast}\Omega
\mid\pi(a)\Omega\,\right\rangle $. Substitution of $L^{\ast}\Omega=-L\Omega=0$
into this identity yields identity (\ref{eq10.9}) for the case $n=1$.

Let $T_{n}$ denote the trace-vector for the trace representative $\tau_{n}$ of
$M_{n}$. Then $\pi_{\omega_{n}}=\pi\otimes\tau_{n}$, and therefore%
\begin{align*}
\left\langle\, \pi_{\omega_{n}}(a\otimes b)\Omega\otimes T_{n}\mid\Omega\otimes
T_{n}\,\right\rangle  &  =\left\langle\, \pi(a)\Omega\otimes\tau_{n}(b)T_{n}%
\mid\Omega\otimes T_{n}\,\right\rangle \\
&  =\left\langle\, \pi(a)\Omega\mid\Omega\,\right\rangle \left\langle\, \tau
_{n}(b)T_{n}\mid T_{n}\,\right\rangle =\omega(a)\operatorname*{tr}%
\nolimits_{n}(b)\\
&  =\omega\otimes\operatorname*{tr}\nolimits_{n}(a\otimes b)=\omega
_{n}(a\otimes b)
\end{align*}
for all $a\in\mathfrak{A}$ and $b\in M_{n}$. Hence $\Omega_{n}=\Omega\otimes
T_{n}$. If we can show that a simple tensor operator $L_{n}$ implements
$\delta_{n}$ in $\pi_{\omega_{n}}$, then identity (\ref{eq10.9}), for
arbitrary $n$, follows from the case $n=1$ which was established above.

However, it is easy to see that the operator $L_{n}=L\otimes I_{n}$ satisfies
the requirements which were listed above. Indeed
\begin{align*}
\pi_{\omega_{n}}(\delta_{n}(a\otimes b))  &  =\pi_{\omega}(\delta
(a))\otimes\tau_{n}(b)\\
&  =(L\pi_{\omega}(a)+\pi_{\omega}(a)L^{\ast})\otimes\tau_{n}(b)\\
&  =L_{n}\pi_{\omega}(a)\otimes\tau_{n}(b)+\pi_{\omega}(a)\otimes\tau
_{n}(b)L_{n}^{\ast}\\
&  =L_{n}\pi_{\omega_{n}}(a\otimes b)+\pi_{\omega_{n}}(a\otimes b)L_{n}^{\ast}%
\end{align*}
for all $a\in\mathfrak{A} $ and $b\in M_{n}$. It follows that $L_{n}$
implements $\delta_{n}$ in $\pi_{\omega_{n}}$.

Only the verification of (\ref{eq10.7}) for $n>1$ then remains. Let $a_{ij}%
\in\mathfrak{A} $ be the matrix entries in some $a\in\mathfrak{A}
_{n}=\mathfrak{A} \otimes M_{n}$. Then the $(i,j)$'th entry $c_{ij}$ in
$a^{\ast}\delta_{n}(a)$ is $\sum_{k=1}^{n}a_{ki}^{\ast}\delta(a_{kj})$. Hence%
\begin{align*}
\omega_{n}(a^{\ast}\delta_{n}(a))  &  =(\omega\otimes\operatorname*{tr}%
\nolimits_{n})(c_{ij})\\
&  =\sum_{i=1}^{n}\omega(c_{ii})=\sum_{i}\sum_{k}\omega(a_{ki}^{\ast}%
\delta(a_{ki})).
\end{align*}
Since $\operatorname{Re}\omega(a_{ki}^{\ast}\delta(a_{ki}))\leq0$,
(\ref{eq10.7}) follows.
\end{proof}

\begin{remark}
\label{Remark8}In the foundations of irreversible statistical thermodynamics
\cite{10,17,18,20,22}, the most conclusive results have been obtained for
dynamical semigroups which are described mathematically as strongly
continuous, completely positive, contraction semigroups $T_{t}$ on the Banach
space $\mathcal{T}(\mathcal{H})$ of all trace-class operators on a given
separable $\infty$-dimensional Hilbert space $\mathcal{H}$. Lindblad \cite{22}
found a formula for the infinitesimal generator
\[
W=\frac{d\,}{dt}T_{t}\bigg|   _{t=0}%
\]
in the case of norm-continuous semigroups, and Davies \cite{11} extended the
results to strongly continuous $T_{t}$ \textup{(}i.e., unbounded generator
$W$\/\textup{)}, satisfying certain side conditions. The condition of
relevance to our paper is the invariance assumption of \cite{11} that
\[
T_{t}^{\prime}(\mathcal{C}(\mathcal{H}))\subset\mathcal{C}(\mathcal{H})
\]
for all $t\in\lbrack0,\infty)$, where $\mathcal{C}(\mathcal{H})$ denotes the
compact operators, and $T_{t}^{\prime}$ the conjugate semigroup on
$B(\mathcal{H})$. Our Theorem \textup{\ref{Theorem1}} does not apply to the
algebra $\mathfrak{A}=B(\mathcal{H})$ since $B(\mathcal{H})^{\prime\prime}$ is
known not to be injective \cite{7}. \textup{(}Of course, $B(\mathcal{H})$ is
injective by Arveson's theorem.\/\textup{)}

However, Theorem \textup{\ref{Theorem1}} combined with the above results
suggests that a $W^{\ast}$-algebra, properly smaller than $B(\mathcal{H})$, is
suitable for quantum dynamics. On the one hand, $B(\mathcal{H})$ \textup{(}or
$\mathcal{T}(\mathcal{H})$ in the conjugate \textup{(}dual\/\textup{)}
formulation\/\textup{)} is too big to accomodate the extensions; and, on the
other hand, the requirement that $\mathcal{C}(\mathcal{H})$ contain the domain
of the generator also appears to be too restrictive.
\end{remark}

\section{Unbounded *-derivations\label{Unb}}

Let $\mathfrak{A} $ be a unital $C^{\ast}$-algebra, and let $\mathcal{D}%
(\delta)$ be a dense $\ast$-subalgebra containing the identity $\openone  $. A
linear transformation $\delta\colon\mathcal{D}(\delta)\rightarrow\mathfrak{A}
$ is said to be a (unbounded) $\ast$-derivation if $\delta(ab)=\delta
(a)b+a\delta(b)$ for $a,b\in\mathcal{D}(\delta)$, and $\delta(a^{\ast}%
)=\delta(a)^{\ast}$ for $a\in\mathcal{D}(\delta)$.

Since, for $\ast$-derivations, one is primarily interested in extensions which
are also $\ast$-derivations, it is natural to work with a two-sided condition
in place of the dissipative notions which were studied in the previous
sections for more general operators. The following such two-sided condition
was suggested by Sakai \cite{29}, and adopted by several authors in subsequent
research on unbounded $\ast$-derivations.

\begin{definition}
\label{Definition9}A $\ast$-derivation $\delta\colon\mathcal{D}(\delta
)\rightarrow\mathfrak{A} $ is said to be well behaved if for all positive
$a\in\mathcal{D}(\delta)$ there is a state $\phi$ on $\mathfrak{A} $ such that
$\phi(a)=\| a\| $ and $\phi(\delta(a))=0$.
\end{definition}

The argument in the previous section yields:

\begin{proposition}
\label{Proposition10}Let $\delta\colon\mathcal{D}(\delta)\rightarrow
\mathfrak{A} $ be a $\ast$-derivation. Then the following four conditions are equivalent:

\begin{enumerate}
\item \label{Proposition10(1)}$\delta$ is well behaved.

\item \label{Proposition10(2)}For all positive $a\in\mathcal{D}(\delta)$, and
for all states $\phi$ on $\mathfrak{A} $ satisfying $\phi(a)=\| a\| $, we have
$\phi(\delta(a))=0$.

\item \label{Proposition10(3)}Each of the operators $\pm\delta$ is dissipative.

\item \label{Proposition10(4)}$\| a+\alpha\delta(a)\| \geq\| a\| $ for all
$\alpha\in\mathbb{R}$ and all $a\in\mathcal{D}(\delta)$.
\end{enumerate}
\end{proposition}

\begin{definition}
\label{Definition11}A $\ast$-derivation $\delta\colon\mathcal{D}%
(\delta)\rightarrow\mathfrak{A} $ is said to be well behaved in the
\emph{matricial sense} if, for each $n=1,2,\dots$, the $\ast$-derivation
$\delta_{n}=\delta\otimes\operatorname*{id}_{n}\colon\mathcal{D}%
(\delta)\otimes M_{n}\rightarrow\mathfrak{A} \otimes M_{n}$ is well behaved.
Recall that $\delta_{n}$ may be regarded as a transformation on $n$-by-$n$
matrices with entries in $\mathfrak{A} $. For such a matrix $a=(a_{ij})$,
$i,j=1,\dots,n$, we have $\delta_{n}(a)=(\delta(a_{ij}))$.
\end{definition}

\begin{theorem}
\label{Theorem12}Every well-behaved $\ast$-derivation is also well behaved in
the matricial sense \textup{(}i.e., completely well behaved\/\textup{)}.
\end{theorem}

\begin{lemma}
\label{Lemma13}Let $\delta\colon\mathcal{D}(\delta)\rightarrow\mathfrak{A} $
be a well-behaved $\ast$-derivation, and let $a\in\mathcal{D}(\delta)$ be
positive. Then there is a state $\phi$ on $\mathfrak{A} $ such that
$\phi(a)=\| a\| $, and $\phi(\delta(b))=0$ for a dense set of elements $b\in
C^{\ast}(a)\cap\mathcal{D}(\delta)$. \textup{(}Here $C^{\ast}(a)$ denotes the
abelian $C^{\ast}$-subalgebra generated by $a$; and every element in $C^{\ast
}(a)$ can be approximated in norm by a sequence of elements $b$ satisfying the
conclusion of the lemma.\/\textup{)}
\end{lemma}

\begin{proof}
[Proofs]The implication (\ref{Proposition10(1)}) $\Rightarrow$
(\ref{Proposition10(2)}) in Proposition \ref{Proposition10} is the key to the
proof of Lemma \ref{Lemma13}. Since functional calculus is also applied, we
shall assume in fact that $\delta$ is closed. By a result of Kishimoto-Sakai
\cite{29} this is no loss of generality. Let $a$ be a positive element in
$\mathcal{D}(\delta)$. Note that the Gelfand-transform sets up an isomorphism
between the $C^{\ast}$-algebras $C^{\ast}(a)$ and $C(\operatorname*{sp}(a))$,
continuous functions on the spectrum of $a$. Let $\lambda_{0}=\mathrm{l.u.b.}%
\operatorname*{sp}(a)$. Then the state $c\rightarrow c(\lambda_{0})$ on
$C(\operatorname*{sp}(a))$ corresponds to a state on $C^{\ast}(a)$ via the
Gelfand-transform. The latter state is then extended to $\mathfrak{A} $ by
Krein's theorem, and the extended state is denoted by $\phi$. It has the
multiplicative property: $\phi(b_{1}b_{2})=\phi(b_{1})\phi(b_{2})$ for
$b_{1},b_{2}\in C^{\ast}(a)$.

Now let $g$ be a non-decreasing (monotone) continuous real function defined on
$\operatorname*{sp}(a)$. Then the Gelfand-transform of $g(a)$ achieves its
maximum at the point $\lambda_{0}$ since the transform of $a$ does. But it is
known that if $g$ is also of class $C^{2}$ (two continuous derivatives) then
$g(a)\in\mathcal{D}(\delta)\cap C^{\ast}(a)$. Hence $\phi(g(a))=\| g(a)\| $.
An application of Proposition \ref{Proposition10}, (\ref{Proposition10(1)})
$\Rightarrow$ (\ref{Proposition10(2)}), then yields the conclusion
\[
\phi(\delta(g(a)))=0.
\]
The restriction of an arbitrary monomial $\lambda^{n}$ to $\operatorname*{sp}%
(a)$ satisfies the conditions listed for $g$. Hence, by Stone-Weierstrass
there is a dense set of elements $b\in C^{\ast}(a)\cap\mathcal{D}(\delta)$
satisfying the conclusion of the lemma. (Alternatively, every positive
function $f$ in $C^{4}$ may be written in the form $f=g_{1}-g_{2}$, with
$g_{1}$ and $g_{2}$ both having the properties listed above for $g$, we
conclude that $\phi(\delta(f(a)))=\phi(\delta(g_{1}(a)))-\phi(\delta
(g_{2}(a)))=0$.)

Now, for each fixed element $a\in\mathcal{D}(\delta)_{+}$ we choose a state
$\phi=\phi_{a}$ and a dense $\ast$-subalgebra $\mathfrak{B} =\mathfrak{B}
_{a}$ of $C^{\ast}(a)$ according to Lemma \ref{Lemma13}; i.e., we require that
$\phi_{a}(\delta(b))=0$ for $b\in\mathfrak{B} _{a}$, as well as $\phi
_{a}(a)=\Vert a\Vert$. Consider the GNS representation of the algebra
$\mathfrak{B} $, resp.\ $\mathfrak{A}$, with representation space
$\mathcal{H}_{\phi}$, resp.\ $\mathcal{K}_{\phi}$, and define:%
\begin{equation}
\mathcal{H}=\sum\nolimits^{\otimes}\mathcal{H}_{\phi},\text{\quad resp.,\quad
}\mathcal{K}=\sum\nolimits^{\otimes}\mathcal{K}_{\phi}. \label{eq11.2}%
\end{equation}
Then $\mathcal{H}$ is a closed subspace of the Hilbert space $\mathcal{K}$,
and we can then define an operator $S$ with dense domain from $\mathcal{H}$ to
$\mathcal{K}$ as follows:%
\begin{equation}
S\pi_{\phi}(b)\Omega_{\phi}=\pi_{\phi}(\delta(b))\Omega_{\phi}\text{\qquad for
}b\in\mathfrak{B} _{\phi}. \label{eq11.3}%
\end{equation}
For vectors $\xi_{1}$ and $\xi_{2}$ in the domain of $S$ we have%
\begin{equation}
\left\langle\, S\xi_{1}\mid\xi_{2}\,\right\rangle +\left\langle\, \xi_{1}\mid
S\xi_{2}\,\right\rangle =0. \label{eq11.4}%
\end{equation}

The verification of (\ref{eq11.4}) may be based on the direct-sum
decomposition (\ref{eq11.2}) above. If $\xi_{i}=\sum_{\phi}^{\otimes}\pi
(b_{i})\Omega_{\phi}\ $for $i=1,2$ and $b_{i}\in\mathfrak{B} _{\phi}$, then
identity (\ref{eq11.4}) reduces to%
\[
\sum\left\langle\, \pi_{\phi}(\delta(b_{1}))\Omega_{\phi}\mid\pi_{\phi}%
(b_{2})\Omega_{\phi}\,\right\rangle +\sum\left\langle\, \pi_{\phi}(b_{1}%
)\Omega_{\phi}\mid\pi_{\phi}(\delta(b_{2}))\Omega_{\phi}\,\right\rangle =0.
\]

The individual terms work out to be:%
\[
\phi(b_{2}^{\ast}\delta(b_{1}))+\phi(\delta(b_{2})^{\ast}b_{1})=\phi
(\delta(b_{2}^{\ast}b_{1}))=0.
\]
Hence, the symmetry condition (\ref{eq11.4}) is hereby reduced to the
conclusion of Lemma \ref{Lemma13} for a given well-behaved derivation $\delta$.

If $P$ denotes the orthogonal projection in $\mathcal{K}$ with range
$\mathcal{H}$, identity (\ref{eq11.4}) implies that the operator
$\xi\rightarrow PS\xi$ may in fact be regarded as a skew symmetric operator in
the Hilbert space $\mathcal{H}$, with dense domain there. We shall also denote
this operator by $S$. The verification of the identity%
\[
\pi(\delta(b))=S\pi(b)-\pi(b)S
\]
is left to the reader.

Following the idea of \S \ref{Con}, we now consider the $\ast$-derivations
$\delta_{n}=\delta\otimes\operatorname*{id}_{n}$ (for each $n=1,2,\dots$)
introduced in Definition \ref{Definition11}. For a given $\ast$-algebra
$\mathfrak{C} $ we denote by $\mathfrak{C} _{n}$ the $\ast$-algebra
$\mathfrak{C} \otimes M_{n}$. Correspondingly, $\ast$-algebras $\mathcal{D}%
(\delta)_{n}$, $\mathfrak{A}_{n}$, and $\mathfrak{B} _{n}$ are defined for
each $n$. Application of the GNS representation to each $\phi_{n}=\phi
\otimes\operatorname*{tr}_{n}$ yields sequences of Hilbert spaces%
\[
\mathcal{H}^{(n)}\subset\mathcal{K}^{(n)}%
\]
as in (\ref{eq11.2}) with each $\mathcal{H}^{(n)}$, resp., $\mathcal{K}^{(n)}%
$, a direct sum of GNS representation spaces associated to $\phi_{n}$.

The calculations in \S \ref{Con} show that the operator $S_{n}=S\otimes I_{n}$
satisfies the $n$'th-order version of (\ref{eq11.3}), that is, (\ref{eq11.3})
holds with the quadruple $S,\pi,\mathfrak{B} ,\delta$ replaced by $S_{n}%
,\pi_{n},\mathfrak{B} _{n},\delta_{n}$. Similarly $\left\langle\,
\smash{S_{n}\xi
_{1}^{(n)}} \bigm|     \smash{\xi_{2}^{(n)}} \,\right\rangle +\left\langle\,
\smash{\xi_{1}^{(n)}} \bigm|    \smash{S_{n}\xi_{2}^{(n)}} \,\right\rangle =0$
for vectors $\xi_{i}^{(n)}$, $i=1,2$, in the respective domains.

Hence Theorem \ref{Theorem6} in \S \ref{Con} implies that each of the
operators $\pm\delta_{n}$ for $n=1,2,\dots$ is dissipative. By Proposition
\ref{Proposition10}, (\ref{Proposition10(3)}) $\Rightarrow$
(\ref{Proposition10(1)}), it follows that $\delta_{n}$ is well behaved,
concluding the proof of Theorem \ref{Theorem12}.
\end{proof}

As an application of the theorem we get the following existence result for
generator extensions of well-behaved $\ast$-derivations $\delta\colon
\mathcal{D}(\delta)\rightarrow\mathfrak{A} $ in \emph{nuclear} $C^{\ast}%
$-algebras $\mathfrak{A} $. Indeed, if $\delta$ is such a $\ast$-derivation,
each of the operators $\pm\delta$ is completely dissipative. Hence, by Theorem
\ref{Theorem1}, there are extensions $\tilde{\delta}_{\pm}\supset\pm\delta$ to
infinitesimal generators of dynamical semigroups $a_{t}^{(\pm)}$ in the
enveloping $W^{\ast}$-algebra $\mathfrak{A} ^{\prime\prime}$.

\begin{acknowledgments}
We thank William Arveson for a helpful email exchange concerning the
existence problem for completely positive semigroups. The book \cite{Arv03a}
offers a different approach.
We are grateful to Brian Treadway for putting a first draft with rewrites and
revisions into beautiful \TeX  .
\end{acknowledgments}


\begin{thebibliography}{99}                                                                                               %


\bibitem {1}Ahiezer, N.I., and Krein, M., \textit{Some Questions in the Theory
of Moments}, Transl. Math. Monographs, Vol. 2, American Mathematical Society,
Providence, 1962.

\bibitem {2}Arveson, W., \textit{Subalgebras of} $C^{\ast}$\textit{-algebras},
Acta Math. \textbf{123} (1969), 141--224.

\bibitem {Arv97}Arveson, W., \textit{Dynamical invariants for noncommutative
flows}, Operator Algebras and Quantum Field Theory: Proceedings of the
conference held at Accademia Nazionale dei Lincei, Roma, Italy, July 1--6,
1996 (S.~Doplicher, R.~Longo, J.E. Roberts, and L.~Zsido, eds.), International
Press, Cambridge, MA, 1997, pp.~476--514.

\bibitem {Arv02a}Arveson, W., \textit{The domain algebra of a CP-semigroup},
Pacific J. Math. \textbf{203} (2002), 67--77.

\bibitem {Arv02b}Arveson, W., \textit{The heat flow of the CCR algebra}, Bull.
London Math. Soc. \textbf{34} (2002), 73--83.

\bibitem {Arv03a}Arveson, W., \textit{Noncommutative Dynamics and {$E$}-Semigroups},
Springer Monographs in Mathematics, Springer-Verlag, New York, 2003.

\bibitem {3}Batty, C.J.K., \textit{Dissipative mappings with approximately
invariant subspaces}, J. Funct. Anal. \textbf{32} (1979), 336--341.

\bibitem {Bha01}Bhat, B.V.R., \textit{Cocycles of CCR Flows}, Mem. Amer. Math.
Soc. \textbf{149} (2001), no.~709.

\bibitem {BrRoII}Bratteli, O., and Robinson, D.W., \textit{{O}perator
{A}lgebras and {Q}uantum {S}tatistical {M}echanics}, 2nd ed., vol.~II,
Springer-Verlag, Berlin--New York, 1996.

\bibitem {4}Bratteli, O., and Kishimoto, A., \textit{Generation of semigroups,
and two-dimensional quantum lattice systems}, J.~Funct. Anal. \textbf{35}
(1980), 344--368.

\bibitem {7}Choi, M.-D., and Effros, E., \textit{Injectivity and operator
spaces}, J.~Funct. Anal. \textbf{24} (1977), 156--209.

\bibitem {8}Choi, M.-D., and Effros, E., \textit{Nuclear }$C^{\ast}%
$\textit{-algebras and injectivity: the general case}, Indiana Univ. Math.~J.
\textbf{26} (1977), 443--446.

\bibitem {9}Connes, A., \textit{Classification of injective factors: Cases}
$\mathrm{II}_{1}$, $\mathrm{II}_{\infty}$, $\mathrm{III}_{\lambda}$,
$\lambda\neq1$, Ann. of Math. (2) \textbf{104} (1976), 73--115.

\bibitem {10}Davies, E.B., \textit{Quantum Theory of Open Systems}, Academic
Press, London--New York, 1976.

\bibitem {11}Davies, E.B., \textit{Generators of dynamical semigroups},
J.~Funct. Anal. \textbf{34} (1979), 421--432.

\bibitem {12}de Leeuw, K., \textit{On the adjoint semigroup and some problems
in the theory of approximation}. Math.~Z. \textbf{73} (1960), 219--234.

\bibitem {13}Dubin, D.A., and Sewell, G.L., \textit{Time translations in the
algebraic formulation of statistical mechanics}, J.~Math. Phys. \textbf{11}
(1970), 2990-2998.

\bibitem {14}Emch, G.G., \textit{Algebraic Methods in Statistical Mechanics
and Quantum Field Theory}, Wiley-Interscience, New York, 1972.

\bibitem {16}Evans, D.E., and Hanche-Olsen, H., \textit{The generators of
positive semigroups}, J.~Funct. Anal. \textbf{32} (1979), 207-212.

\bibitem {17}Evans, D.E., and Lewis, J.T., \textit{Dilations of Irreversible
Evolutions in Algebraic Quantum Theory}, Comm. Dublin Inst. Adv. Studies Ser.
A No. 24 (1977).

\bibitem {18}Evans, D.E., \textit{A review on semigroups of completely
positive maps}, Mathematical Problems in Theoretical Physics (Lausanne, 1979)
(K.~Osterwalder, ed.), Lecture Notes in Physics, vol.~116, Springer-Verlag,
Berlin--Heidelberg--New York, 1980, pp.~400--406.

\bibitem {FMRR00}Ferrari, P.A., Maes, C., Ramos, L., and Redig, F., \textit{On
the hydrodynamic equilibrium of a rod in a lattice fluid}, J. Phys. A
\textbf{33} (2000), 4725--4740.

\bibitem {19}Hille, E., and Phillips, R.S., \textit{Functional Analysis and
Semi-groups}, American Mathematical Society, Providence, 1957.

\bibitem {20}Ingarden, R.S., and Kossakowski, A., \textit{On the connection of
nonequilibrium information theormodynamics with non-Hamiltonian quantum
mechanics of open systems}, Ann. Physics \textbf{89} (1975), 451--485.

\bibitem {JaPi01}Jak{\v{s}}i{\'{c}, V.,} and Pillet, C.-A., \textit{On entropy
production in quantum statistical mechanics}, Comm. Math. Phys. \textbf{217}
(2001), 285--293.

\bibitem {21}Jorgensen, P.E.T., \textit{Approximately invariant subspaces for
unbounded linear operators}, II, Math. Ann. \textbf{227} (1977), 177--182.

\bibitem {Lin00}Lindblad, G{.}, \textit{Cloning the quantum oscillator}, J.
Phys. A \textbf{33} (2000), 5059--5076.

\bibitem {22}Lindblad, G., \textit{On the generators of quantum dynamical
semigroups}, Comm. Math. Phys. \textbf{48} (1976), 119--130.

\bibitem {23}Lumer, G., and Phillips, R.S., \textit{Dissipative operators in a
Banach space}, Pacific J. Math. \textbf{11} (1961), 679--698.

\bibitem {24}Phillips, R.S., \textit{Dissipative operators and hyperbolic
systems of partial differential equations}, Trans. Amer. Math. Soc.
\textbf{90} (1959), 193--254.

\bibitem {25}Powers, R.T., and Sakai, S., \textit{Unbounded derivations in
operator algebras}, J.~Funct. Anal. \textbf{19} (1975), 81--95.

\bibitem {26}Pulvinenti, M., and Tirozzi, B., \textit{Time evolution of a
quantum lattice system}, Comm. Math. Phys. \textbf{30} (1973), 83--98.

\bibitem {27}Ruelle, D., \textit{Statistical Mechanics: Rigorous Results}, W.
A. Benjamin, New York--Amsterdam, 1969.

\bibitem {Rue01}Ruelle, D., \textit{Topics in quantum statistical mechanics
and operator algebras}, preprint, 2001.

\bibitem {28}Sakai, S., $C^{\ast}$\textit{-Algebras and} $W^{\ast}%
$\textit{-Algebras}, Ergebnisse der Mathematik und ihrer Grenzgebiete, Band
60, Springer-Verlag, Berlin--Heidelberg--New York, 1971.

\bibitem {29}Sakai, S., \textit{Recent developments in the theory of unbounded
derivations in} $C^{\ast}$\textit{-algebras},
$C^{\ast}$-Algebras and
Applications to Physics (Los Angeles, 1977) (H.~Araki and R.V. Kadison, eds.),
Lecture Notes in Mathematics, vol.~650, Springer-Verlag,
Berlin--Heidelberg--New York, 1978, pp.~85--122.

\bibitem {SW00}Salmhofer, M., and Wieczerkowski, C., \textit{Positivity and
convergence in fermionic quantum field theory}, J. Statist. Phys. \textbf{99},
(2000), 557--586.

\bibitem {30}Tomiyama, J., \textit{On the projection of norm one in} $W^{\ast
}$\textit{-algebras}, Proc. Japan Acad. \textbf{33} (1957), 608--612.
\end{thebibliography}
\end{document}